\documentclass{amsproc}
\usepackage{amsmath,amssymb,amsfonts}
\usepackage[mathscr]{eucal}

\usepackage{enumerate}

\theoremstyle{plain}
\newtheorem{theorem}{Theorem}[section]

\theoremstyle{definition}
\newtheorem{definition}[theorem]{Definition}
\newtheorem{example}[theorem]{Example}

\numberwithin{equation}{section}

\input xy
\xyoption{all}

\newcommand{\Deltaop}{{\bf \Delta}^{op}}

\newcommand{\hocolim}{\text{hocolim}}

\newcommand{\nerve}{\text{nerve}}

\newcommand{\we}{\text{we}}

\newcommand{\Aut}{\text{Aut}}

\newcommand{\SSets}{\mathcal{SS}ets}

\newcommand{\Sets}{\mathcal Sets}
\newcommand{\map}{\text{map}}
\newcommand{\ob}{\text{ob}}
\newcommand{\hoequiv}{\text{hoequiv}}
\newcommand{\css}{\mathcal{CSS}}

\begin{document}

\title[Homotopy colimits]{Homotopy colimits of model categories}

\author[J.E. Bergner]{Julia E. Bergner}
\address{Department of Mathematics, University of California, Riverside, CA 92521}

\email{bergnerj@member.ams.org}
\thanks{The author was partially supported by NSF grant DMS-1105766.}


\date{\today}

\begin{abstract}
Building on a previous definition of homotopy limit of model categories, we give a definition of homotopy colimit of model categories.  Using the complete Segal space model for homotopy theories, we verify that this definition corresponds to the model-category-theoretic definition in that setting.
\end{abstract}

\maketitle

\section{Introduction}

Model categories, first developed by Quillen in \cite{quillen}, have been an important tool in extending the ideas of homotopy theory from the world of topological spaces to other areas of mathematics.  While the basic data of a homotopy theory consists of a category together with some choice of weak equivalences, having the additional structure of a model category allows one to make homotopy-invariant constructions precise.  For example, one can define homotopy limits and colimits in a model category.

A more recent perspective on homotopy theory, begun by Dwyer and Kan in their work on simplicial localizations \cite{dkfncxes}, \cite{dksimploc}, takes homotopy theories themselves as the objects of study.  Therefore, one could try to make the same kinds of structures that have been investigated \emph{within} model categories applicable to some category \emph{of} model categories.  For example, what should a homotopy limit or colimit of a diagram of model categories be?

Answering such questions has been much easier in the context of more general homotopy theories, also known as $(\infty,1)$-categories.  Categories with weak equivalences, simplicial categories, quasi-categories, Segal categories, and complete Segal spaces are all ways of modeling homotopy theories as mathematical objects; there are model categories for each of these different models which are all Quillen equivalent \cite{bk}, \cite{simpcat}, \cite{thesis}, \cite{joyal}, \cite{jt}, \cite{lurie}, \cite{pell}.  Therefore, in any one of these settings, one can use standard model category techniques to understand what is meant by a homotopy limit or homotopy colimit of homotopy theories.

However, it is worth trying to make these constructions in the more restrictive world where the objects are actually model categories.  There is no known model category of model categories, so there is no immediate approach to take.  A definition of homotopy fiber product of model categories was given in \cite{toendha}, however.  In \cite{fiberprod} we used an explicit functor from model categories to complete Segal spaces to verify that this definition did in fact correspond to a homotopy pullback of homotopy theories in that setting.  In \cite{holim}, we extended this definition, and the proof of its validity, to more general homotopy limits of model categories.

Here, we consider homotopy colimits of model categories.  The situation is somewhat worse, in that there is very little expectation that we actually get a model structure.  For example, just taking a coproduct, or disjoint union of model categories, does not result in a model category; we certainly don't have products and coproducts unless all objects come from the same original model category.

As with homotopy limits of model categories, it is expected that many applications fit into this framework.  In current work with Robertson and Salch, we use this construction to understand a topological triangulated orbit category as a homotopy coequalizer of certain stable model categories \cite{tto}.

In Section 2, we give a review of complete Segal spaces.  We treat the case of homotopy pushouts of model categories in Section 3, then go on to the general case in Section 4, establishing that a homotopy colimit of model categories translates into a genuine homotopy colimit of complete Segal spaces.

\section{Preliminaries on complete Segal spaces}

In this section we give the necessary background on complete Segal spaces for the arguments later in the paper.

Recall that a simplicial set is a functor $\Deltaop \rightarrow \Sets$, where $\Deltaop$ is the simplicial indexing category.  We denote the category of simplicial sets by $\SSets$; it has a model structure Quillen equivalent to the usual model structure on topological spaces.  A \emph{simplicial space} is a functor $\Deltaop \rightarrow \SSets$.  We refer the reader to \cite{gj} for more details about simplicial sets and other simplicial objects.  In particular, the category of simplicial spaces can be given the Reedy model structure, in which the weak equivalences are given by levelwise weak equivalences of simplicial sets, and cofibrations are monomorphisms \cite[15.8.7]{hirsch}.

\begin{definition} \cite[4.1]{rezk}
A \emph{Segal space} is a Reedy fibrant simplicial space $W$ such that the Segal maps
\[ \varphi_n \colon W_n \rightarrow \underbrace{W_1 \times_{W_0} \cdots \times_{W_0} W_1}_n \] are weak equivalences of simplicial sets for all $n \geq 2$.
\end{definition}

Given a Segal space $W$, its ``objects" are defined by $\ob(W)= W_{0,0}$, and, between any two objects $x$ and $y$, the ``mapping space" $\map_W(x,y)$, defined to be the homotopy fiber of the map $W_1 \rightarrow W_0 \times W_0$ given by the two face maps $W_1 \rightarrow W_0$.  The condition on the Segal maps guarantees a notion of $n$-fold composition of mapping spaces, up to homotopy.  Using this composition, we can define ``homotopy equivalences", and then consider of the subspace of $W_1$ whose components contain homotopy equivalences, denoted by $W_{\hoequiv}$.  Then the degeneracy map $s_0 \colon W_0 \rightarrow W_1$ factors through $W_{\hoequiv}$.

\begin{definition} \cite[\S 6]{rezk}
A \emph{complete Segal space} is a Segal space $W$ such that the map $W_0 \rightarrow W_{\hoequiv}$ is a weak equivalence of simplicial sets.
\end{definition}

\begin{theorem} \cite[\S 7]{rezk}
There is a model category structure $\css$ on the category of simplicial spaces, obtained as a localization of the Reedy model structure such that:
\begin{enumerate}
\item the fibrant objects are the complete Segal spaces,

\item all objects are cofibrant, and

\item the weak equivalences between complete Segal spaces are levelwise weak equivalences of simplicial sets.
\end{enumerate}
\end{theorem}

In particular, colimits of complete Segal spaces can be taken as levelwise colimits of simplicial sets.

Rezk defines a functor which we denote $L_C$ from the category of model categories and left Quillen functors to the category of simplicial spaces; given a model category $\mathcal M$, we have that
\[ L_C(\mathcal M)_n = \nerve(\we(\mathcal M^{[n]})). \]  Here, $\mathcal M^{[n]}$ is the category of maps $[n] \rightarrow \mathcal M$, and $\we(\mathcal M^{[n]})$ denotes the subcategory of $\mathcal M^{[n]}$ whose morphisms are the weak equivalences.  While the resulting simplicial space is not in general Reedy fibrant, and hence not a complete Segal space, Rezk proves that taking a Reedy fibrant replacement is sufficient to obtain a complete Segal space \cite[8.3]{rezk}.  For the rest of this paper we assume that the functor $L_C$ includes composition with this Reedy fibrant replacement and therefore assigns a complete Segal space to any model category.  In fact, this construction can be applied to any category with weak equivalences.

A difficulty with this definition is the fact that it is only a well-defined functor on the category whose objects are model categories and whose morphisms preserve weak equivalences, not on the category of model categories with morphisms left Quillen functors.  Instead, consider $\mathcal M^c$, the full subcategory of $\mathcal M$ whose objects are cofibrant.  While $\mathcal M^c$ may no longer have the structure of a model category, it is still a category with weak equivalences.  Thus, we define
\[ L_C(\mathcal M)_n = \nerve(\we((\mathcal M^c)^{[n]})). \]  Each space in this diagram is weakly equivalent to the one given by the previous definition, and now the construction is functorial on the category of model categories with morphisms the left Quillen functors.  To consider right Quillen functors instead, take the full subcategory of fibrant objects, $\mathcal M^f$, rather than $\mathcal M^c$.

It remains to give a description of the image of a model category $\mathcal M$ under $L_C$.  We begin with some notation.  Given a simplicial
monoid $M$, there is a classifying complex of $M$, a simplicial set whose geometric realization is the classifying space $BM$ \cite[V.4.4]{gj}, \cite{may}.  We simply write
$BM$ for the classifying complex of $M$.  We also consider disjoint unions of simplicial monoids; in this case the classifying complex is taken in the category of simplicial categories, rather than in simplicial monoids.

\begin{theorem} \cite[7.3]{css}  \label{baut}
Let $\mathcal M$ be a model category.  For $x$ an object of $\mathcal M$ denote by $\langle x \rangle$ the weak equivalence class of $x$ in $\mathcal M$, and denote by $Aut^h(x)$ the simplicial monoid of self weak equivalences of $x$.  Similarly, let $\langle \alpha \colon x \rightarrow y \rangle$ denote the weak equivalence class of $\alpha$ in $\mathcal M^{[1]}$ and $Aut^h(\alpha)$ its respective simplicial monoid of self weak equivalences.  Up to weak equivalence in the model category $\css$, the complete Segal space $L_C(\mathcal M)$ looks like
\[ \coprod_{\langle x \rangle} BAut^h(x) \Leftarrow \coprod_{\langle \alpha \colon x \rightarrow y \rangle} BAut^h(\alpha) \Lleftarrow \cdots. \]
\end{theorem}

The reference (Theorem 7.3 of \cite{css}) gives a characterization of the complete Segal space arising from a simplicial category, not from a model category.  However, the results of \S 6 of that same paper (specifically, the composite of Theorems 6.2 and 6.4) allow for translating it to the theorem as stated here.  We note additionally that the same characterization applies when $\mathcal M$ is simply a category with weak equivalences, but in this case we cannot assume that $L_C$ is a functor without further assumptions about preservation of weak equivalences.

\section{Working example: Homotopy pushouts}

In this section, we give an informal treatment of the special case of homotopy pushouts in order to get an intuitive sense of how homotopy colimits of model categories should be defined.  We give the formal definition and its justification in the next section.

Consider the diagram of left Quillen functors
\[ \xymatrix{\mathcal M_3 \ar[r]^{F_1} \ar[d]_{F_2} & \mathcal M_1 \ar[d] \\
\mathcal M_2 \ar[r] & \mathcal P. } \]
We expect that the homotopy pushout $\mathcal P$ should be a ``quotient" of $\mathcal M_1 \amalg \mathcal M_2$, where an object $x_1$ of $\mathcal M_1$ is identified with an object $x_2$ of $\mathcal M_2$ if there exists an object $x_3$ in $\mathcal M_3$ together with weak equivalences $F_1(x_3) \rightarrow x_1$ in $\mathcal M_1$ and $F_2(x_3) \rightarrow x_2$ in $\mathcal M_2$.  By ``identify", we mean there should be a string of weak equivalences connecting $x_1$ and $x_2$.

Similarly, if we had a diagram of right adjoint functors
\[ \xymatrix{\mathcal M_3 \ar[r]^{G_1} \ar[d]_{G_2} & \mathcal M_1 \ar[d] \\
\mathcal M_2 \ar[r] & \mathcal P } \]
for which $\mathcal P$ is a homotopy pushout, we would expect to identify $x_1$ and $x_2$ if there exists $x_3$ in $\mathcal M_3$ equipped with weak equivalences $x_1 \rightarrow G_1(x_3)$ in $\mathcal M_1$ and $x_2 \rightarrow G_2(x_3)$ in $\mathcal M_2$.

Notice that we have implicitly simplified the situation here, using the fact that $\mathcal M_3$ is an initial object in our diagram.  To be more rigorous, namely to give a precise identification of objects, we can think of taking the disjoint union $\mathcal M_1 \amalg \mathcal M_2 \amalg \mathcal M_3$ and adding in weak equivalences $x_3 \rightarrow x_1$ and $x_3 \rightarrow x_2$ (or with direction reversed in the right adjoint case), together with all generated composites.  Every object of $\mathcal M_3$ is identified both with an object of $\mathcal M_1$ and with an object of $\mathcal M_2$, so we can informally think of the identification as being between $\mathcal M_1$ and $\mathcal M_2$.  We need to include $\mathcal M_3$ formally, and in the case of more general homotopy colimits, we cannot make such an assumption anyway.

\begin{example}
Let $F \colon \mathcal M \rightarrow \mathcal N$ be a left Quillen functor.  Then we can form the homotopy pushout diagram
\[ \xymatrix{\mathcal M \ar[r]^F \ar[d] & \mathcal N \ar[d] \\
\ast \ar[r] & \mathcal P} \] where $\mathcal P$ denotes the category described above, and $\ast$ denotes the category with a single object and identity morphism only.  Then an object of $\mathcal N$ becomes weakly equivalent to the object of $\ast$ if it is weakly equivalent to an object in the image of $F$.  In other words, the objects which are not identified are those not coming from $\mathcal M$, even up to homotopy, so it makes sense to think of $\mathcal P$ as the homotopy cofiber of the functor $F$.
\end{example}

\section{More general homotopy colimits}

In this section we give a formal definition of homotopy colimits of model categories and justify it by comparing to homotopy colimits in the context of complete Segal spaces.

\begin{definition}
Let $\mathcal D$ be a small category, and $\mathcal M$ a $\mathcal D$-shaped diagram of left Quillen functors $F_{\alpha, \beta}^\theta \colon \mathcal M_\alpha \rightarrow \mathcal M_\beta$.  (Here the superscript $\theta$ allows us to distinguish between different arrows $\alpha \rightarrow \beta$ in $\mathcal D$.) Then the \emph{homotopy colimit} of $\mathcal M$, denoted by $\mathcal Colim_\alpha\mathcal M_\alpha$, is defined to be the category obtained from the disjoint union of the model categories in $\mathcal M$ by inserting weak equivalences $x_\beta \rightarrow x_\alpha$ between objects $x_\alpha$ in $\mathcal M_\alpha$ and $x_\beta$ in $\mathcal M_\beta$ if there exists a weak equivalence $F_{\alpha, \beta}^\theta(x_\alpha) \rightarrow x_\beta$ in $\mathcal M_\beta$.  We further assume that, if such a weak equivalence already exists (in the case where $\alpha=\beta$), we do not add an additional one, and that we impose the appropriate relation on composites: if there exist two weak equivalences $F_{\alpha, \beta}^\theta(x_\alpha) \rightarrow x_\beta$ and $F_{\beta, \gamma}^\psi(x_\beta) \rightarrow x_\gamma$, then the two possible ways of obtaining weak equivalences $x_\alpha \rightarrow x_\gamma$ are identified.

If $\mathcal M$ is instead a $\mathcal D$-shaped diagram of right Quillen functors $G_{\alpha, \beta}^\theta \colon \mathcal M_\alpha \rightarrow \mathcal M_\beta$, then a weak equivalence $x_\alpha \rightarrow x_\beta$ is included if there exists a weak equivalence $x_\beta \rightarrow G_{\alpha, \beta}^\theta(x_\alpha)$ in $\mathcal M_\beta$.
\end{definition}

Observe that this definition can be regarded as a Grothendieck-type construction on the category of model categories.

Ideally, we would like the homotopy colimit of a diagram of model categories to be a model category itself.  Even for homotopy limits of model categories, which are subcategories of the product of model categories, fairly strict assumptions are needed to guarantee a model structure \cite{holim}.  Coproducts of model categories are still less likely to be sensible candidates for a model structure, so we cannot really expect homotopy colimits of model categories to be model categories.  However, they still form categories with weak equivalences.  The weak equivalences are defined to be the maps which are weak equivalences in their original categories, together with those which we have included, and composites needed for the 2-out-of-3 property to hold.

To justify our definition of homotopy colimit, we can use the functor $L_C$ defined in the previous section to translate to the model category of complete Segal spaces, where homotopy colimits are defined.  Then we have the following result.

\begin{theorem}
The map $L_C(\mathcal Colim_\alpha \mathcal M_\alpha) \rightarrow \hocolim_\alpha (L_C \mathcal M_\alpha)$ is a weak equivalence in the complete Segal space model structure.
\end{theorem}

In other words, whether we take a homotopy colimit of model categories in the above sense and then translate to a complete Segal space, or translate to a diagram of complete Segal spaces and take the usual homotopy colimit, we get weakly equivalent results.  The idea of the proof is analogous to that given in \cite{holim} for homotopy limits.

\begin{proof}
We first consider the 0-space of both complete Segal spaces.  We have that
\[ \begin{aligned}
\left( \hocolim_\alpha L_C \mathcal M_\alpha \right)_0 & \simeq \hocolim_\alpha \left( \coprod_{\langle x_\alpha \rangle} B\Aut^h(x_\alpha) \right) \\
& \simeq \coprod_{\langle x_\alpha \rangle} \hocolim_\alpha (B\Aut^h(x_\alpha)).
\end{aligned} \]
Using the construction of homotopy colimits as in \cite{bouskan}, this space is given by taking the coproduct
\[ \coprod_\alpha \left( \coprod_{\langle x_\alpha \rangle} B\Aut^h(x_\alpha) \right) \] and identifying equivalences of automorphisms $a_\alpha \colon x_\alpha \rightarrow x_\alpha$ with $a_\beta \colon x_\beta \rightarrow x_\beta$ (or rather, their images in their respective classifying spaces) via gluing 1-simplices if there exists a commutative diagram of the form
\[ \xymatrix{F_{\alpha, \beta}^\theta(x_\beta) \ar[r]^-\simeq \ar[d]_{F_{\alpha, \beta}^\theta(a_\beta)} & x_\alpha \ar[d]^{a_\alpha} \\
F_{\alpha, \beta}^\theta(x_\beta) \ar[r]^-\simeq & x_\alpha} \]
in $\mathcal M_\alpha$, and analogously for higher simplices for compositions.  However, this description is exactly that of the zero space of $L_C \mathcal Colim_\alpha (\mathcal M_\alpha)$.

A similar argument can be made on the level of morphisms to show that the two simplicial spaces have weakly equivalent 1-spaces; the higher-degree spaces are all determined, by the Segal condition.
\end{proof}

\begin{example}
Let $\mathcal T$ be a triangulated category and $F \colon \mathcal T \rightarrow \mathcal T$ a self-equivalence.  In \cite{keller}, Keller considers the orbit category $\mathcal T/F$, for $\mathcal T$ an algebraic triangulated category, and gives conditions under which it still has a triangulated structure.  His primary example is that of the cluster category.  In \cite{tto}, we consider the case where $\mathcal T$ is a topological triangulated category, or the homotopy category of a stable model category or more general cofibration category.  The definition given here allows for a definition of the orbit category associated to a stable model category equipped with a self-equivalence, and in particular a notion of topological cluster category.
\end{example}

\bibliographystyle{amsplain}

\end{document}